\documentclass[12pt,oneside,reqno]{amsart}
\usepackage{mathrsfs}
\usepackage{graphics}
\usepackage{amssymb}
\newcommand{\dif}{\mathrm{d}}

\newcommand{\be}{\begin{eqnarray}}
\newcommand{\ee}{\end{eqnarray}}
\newcommand{\ce}{\begin{eqnarray*}}
\newcommand{\de}{\end{eqnarray*}}
\newtheorem{theorem}{Theorem}[section]
\newtheorem{lemma}[theorem]{Lemma}
\newtheorem{remark}[theorem]{Remark}
\newtheorem{definition}[theorem]{Definition}
\newtheorem{proposition}[theorem]{Proposition}
\newtheorem{example}[theorem]{Example}
\newtheorem{corollary}[theorem]{Corollary}

\newcommand{\PX}{{\Bbb{P}}}

\def\[{{\Big[}}
\def\]{{\Big]}}
\def\<{{\langle}}
\def\>{{\rangle}}
\def\({{\Big(}}
\def\){{\Big)}}

\def\bt{\begin{theorem}}
\def\et{\end{theorem}}
\def\bl{\begin{lemma}}
\def\el{\end{lemma}}
\def\br{\begin{remark}}
\def\er{\end{remark}}
\def\bx{\begin{Examples}}
\def\ex{\end{Examples}}
\def\bd{\begin{definition}}
\def\ed{\end{definition}}
\def\bp{\begin{proposition}}
\def\ep{\end{proposition}}
\def\bc{\begin{corollary}}
\def\ec{\end{corollary}}

\def\cB{{\mathcal B}}

\def\cF{{\mathcal F}}

\def\cL{{\mathcal L}}

\def\cX{{\mathcal X}}

\def\mE{{\mathbb E}}

\def\mN{{\mathbb N}}

\def\mP{{\mathbb P}}

\def\mR{{\mathbb R}}

\def\mX{{\mathbb X}}

\def\geq{\geqslant}
\def\leq{\leqslant}

\begin{document}

\allowdisplaybreaks

\title{A Multiplicative Ergodic Theorem for Discontinuous Random Dynamical Systems and Applications*}

\author{Huijie Qiao$^1$ and Jinqiao Duan$^2$}



\thanks{*This work was partly supported by NSF of China (No. 11001051, 11371352).}


\date{}

\dedicatory{1. Department of Mathematics,
Southeast University\\
Nanjing, Jiangsu 211189,  China\\
hjqiaogean@seu.edu.cn\\
2. Department of Applied Mathematics, Illinois Institute of Technology\\
Chicago, IL 60616, USA\\
duan@iit.edu}

\begin{abstract}
Motivated by studying stochastic systems with non-Gaussian L\'evy noise, spectral properties
for a type of linear cocycles are considered.   These linear cocycles have countable jump
discontinuities in time. A  multiplicative ergodic theorem is proved for such  linear cocycles.
Then, the result is illustrated for two linear stochastic systems with general L\'evy motions.
\end{abstract}

\maketitle \rm

\section{Introduction}

Multiplicative ergodic theorems (METs) provide  a spectral theory for linear cocycles, which are often solution mappings for linear
stochastic differential equations. The type of theorems provides a stochastic counterpart for deterministic linear algebra, with
spectral objects such as invariant subspaces, exponential growth rate or Lyapunov exponents \cite{la}. These spectral objects
establish a foundation for investigating nonlinear stochastic dynamical systems.

METs for linear cocycles have been summarized in \cite{la}, where these linear cocycles are required to be continuous
in time variable $t$. These linear cocycles often come from the solution mappings of linear stochastic differential equations
(SDE) with (Gaussian) Brownian motions; many authors have considered Lyapunov exponents for these equations, such as \cite{epvw, mpvw}.
METs for linear cocycles with respect to time $t$ in infinite dimensional space have recently been proved in \cite{ll} (also see references therein).
Lyapunov exponents for linear stochastic functional differential equations in \cite{seam, ms1, ms2} have been also investigated.

However, METs for two-sided linear cocycles which are discontinuous (right-continuous with left limits) in $t$ are not available in literature. Although Oseledet \cite{ose} deduced a MET for one-sided linear cocycles which were only measurable in time and
satisfied two integrability conditions, these conditions were strong and not easy to justify. And Li-Blankenship \cite{lb} only studied Lyapunov exponents for one-sided linear systems with Poisson noise. Besides, Mohammed and Scheutzow
\cite{ms1} studied a MET for linear stochastic functional differential equations driven by discontinuous semimartingales for $t\in\mR_+$. There,
they required that the martingale parts of these semimartingales were continuous in $t$.

In this paper, we prove a MET for two-sided linear cocycles discontinuous in $t$ under two weaker integrability conditions. These cocycles
have countable jump discontinuities in time. On one side, by the MET, exponential stability of a number of SDEs with non-Gaussian
L\'evy noise is solved. On the other side, METs in \cite{la, ms1} are generalized by the MET. Besides, our proof is different from Oseledet's one in \cite{ose},
but follows \cite{Franco, ll} with the help of Lemma \ref{met1}.

This paper is arranged as follows. In Section \ref{prelim}, we introduce linear cocycles discontinuous in $t$, flags, and
L\'evy processes for $t\geq 0$. Two motivated examples are placed in Section \ref{exam}. A multiplicative ergodic
theorem (Theorem \ref{met2}) for two-sided linear cocycles which are discontinuous in $t$ is proved in Section \ref{met}. Moreover,
we illustrate this MET by applying it to the two examples in Section \ref{exam} .

\medskip

The following convention will be used throughout the paper: $C$ with or without indices will denote different positive constants
(depending on the indices) whose values vary.

\section{Preliminaries}\label{prelim}

In this section, we recall  basic concepts and facts that will be
needed throughout the paper.

In the following, $|\cdot|$ stands for the length of a vector in
$\mR^d$, $\|\cdot\|$ denotes the Hilbert-Schmidt norm of a matrix or
the norm of a linear operator and $\<\cdot,\cdot\>$ is the usual
scalar product in $\mR^d$.

\subsection{Probability space}

Let $D(\mR,\mR^d)$ be the set of all c\`adl\`ag functions $f$ defined on $\mR$ with values
in $\mR^d$ and $f(0)=0$. We take $\Omega :=D(\mR,\mR^d)$, which
will be the canonical sample space for stochastic differential
equations with two-sided L\'evy motions. It can be made a complete
and separable metric space when endowed with the Skorohod metric
$\rho$ as in \cite{hwy}: for   $x,y\in\Omega$, \ce
\rho(x,y):=\inf\limits_{\lambda\in\Lambda}\left\{\sup\limits_{s\neq
t}\left|\log\frac{\lambda(t)-\lambda(s)}{t-s}\right|
+\sum\limits_{m=1}^\infty\frac{1}{2^m}\min\Big\{1,
\rho^\circ(x^m,y^m)\Big\}\right\}, \de where $x^m(t):=g_m(t)x(t)$
and $y^m(t):=g_m(t)y(t)$ with \ce g_m(t):=\left\{\begin{array}{c}
1, ~\quad\mbox{if}~ |t|\leq m,\quad\\
m+1-|t|,\quad\mbox{if}~m<|t|<m+1,\quad\\
0, \qquad \mbox{if}~|t|\geq m+1,
\end{array}
\right. \de and
\ce
\rho^\circ(x,y):=\sup\limits_{t\in\mR}\left|x(t)-y(\lambda(t))\right|.
\de
Here $\Lambda$ denotes the
set of strictly increasing and continuous functions $\lambda$ from $\mR$ to
$\mR$ with $\lambda(0)=0$. We identify a function $\omega(t)$  with a (canonical) sample
$\omega$ in the sample space $\Omega$.

The Borel $\sigma$-field  in the sample space $\Omega$,  under the
topology induced by the Skorohod metric $\rho$, is denoted by $\cF$.
Note that $\cF=\sigma(\omega(t),t\in\mR)$, as known in
 \cite{hwy}. Let $\mP$ be the unique probability measure which makes the
canonical process a L\'evy process for $t\in\mR$ (Definition
\ref{levy1}). And we have the complete natural
filtration $\cF_s^t:=\sigma(\omega(u): s\leq u\leq
t)\vee\mathcal{N}$ for $s\leq t$ with respect to $\mP$. Here
$\mathcal{N}$ is the set of all null events under $\PX$.

\subsection{Definition of random dynamical systems (RDSs), discontinuous RDSs and linear RDSs}

Define for each $t\in\mR$ \ce
(\theta_t\omega)(\cdot)=\omega(t+\cdot)-\omega(t), \quad
\omega\in\Omega. \de Then $\{\theta_t\}$ is a one-parameter group
(or a \emph{flow}, or a \emph{deterministic dynamical system}) on
$\Omega$. In fact, $\Omega$ is invariant with respect to
$\{\theta_t\}$, i.e.
$$
\theta_t^{-1}\Omega=\Omega, \quad ~\mbox{for~all}~ t\in\mR,
$$
and $\mP$ is $\{\theta_t\}$-invariant, i.e.
$$
\mP(\theta_t^{-1}(B))=\mP(B), \qquad ~\mbox{for~ all}~ B\in\cF,
t\in\mR.
$$
Thus $(\Omega,\cF,\mP,(\theta_t)_{t\in\mR})$ is a metric dynamical
system (DS), or also called a driving dynamical system. The metric
DS $(\Omega,\cF,\mP,(\theta_t)_{t\in\mR})$ is called ergodic, if all
measurable $\{\theta_t\}$-invariant sets have probability $0$ or
$1$. (see \cite{la})

\bd Let $(\mX,\mathcal{X})$ be a measurable space. For a mapping \ce
&&\varphi: \mR\times\Omega\times\mX\mapsto\mX, \quad
(t,\omega,x)\mapsto\varphi(t,\omega,x),\\
&&\varphi(t,\omega):=\varphi(t,\omega,\cdot): \mX\mapsto\mX, \de

(i) Measurability: $\varphi$ is $\cB(\mR)\bigotimes\cF\bigotimes\cX/\cX$-measurable;

(ii) Cocycle property: $\varphi(t,\omega)$ satisfies the following conditions
\be
\varphi(0,\omega)&=&id_{\mX},  \label{perfect coc1}\\
\varphi(t+s,\omega)&=&\varphi(t,\theta_s\omega)\circ\varphi(s,\omega),
\label{perfect coc2} \ee for all $s,t\in\mR$ and $\omega\in\Omega$;
$\varphi$ is called as a random dynamical system (RDS). Sometimes we simply call RDS $\varphi$ a cocycle.
\ed

\bd
A RDS $\varphi$ is called as the discontinuous RDS if it is c\`adl\`ag (right-continuous with left limits) in $t$.
\ed

\bd
A RDS $\varphi$ is called as the linear RDS if
for each $t\in\mR$ and $\omega\in\Omega$, $\varphi(t,\omega)$ is a
linear operator in $\mX$. Sometimes we simply call linear RDS $\varphi$ a linear cocycle.
\ed

\subsection{Flags and related metrics}\label{fm}

Let us introduce the definition for a flag of type $\tau$. Let
$\tau$ be a $p$-dimensional vector with positive integer components
such that $\tau=(d_p,\cdots,d_1)$ and $1\leq d_p<\cdots<d_1=d$. A
flag of type $\tau$ in $\mR^d$ is a sequence of subspaces
$F=(V_p,\cdots,V_1)$ such that $V_p\subset\cdots\subset V_1=\mR^d$
and $\dim V_i=d_i$ for all $i$. The set for all flags of type $\tau$
constitutes the space of flags $F_\tau(d)$. Moreover, $F_\tau(d)$
can be given a structure of a compact $C^\infty$ manifold in a
natural way (\cite{hu}). And the flag manifold $F_\tau(d)$ can be
endowed with a complete metric $\delta$ as follows (\cite{gm}):

Let $U_p$ be equal to $V_p$ and $U_i$ be the orthogonal complement
of $V_{i+1}$ in $V_i$, $i=p-1,\cdots,1$, so that \ce
V_i=U_p\oplus\dots\oplus U_i, \quad i=p,\dots,1. \de Define for any
$F=(V_p,\cdots,V_1), \tilde{F}=(\tilde{V}_p,\cdots,\tilde{V}_1)\in
F_\tau(d)$ \ce \tilde{\rho}(F,\tilde{F}):=\max_{i\neq j, x\in U_i,
y\in\tilde{U}_j;\atop
\|x\|=\|y\|=1}|\<x,y\>|^{h/|\lambda_i-\lambda_j|}, \de where
$\lambda_1, \lambda_2, \cdots, \lambda_p$ and $h$ are real numbers
which satisfy $\lambda_i\neq\lambda_j$ for $i\neq j$ and \ce
h^{-1}|\lambda_i-\lambda_j|\geq d-1 ~\mbox{for}~ i\neq j. \de By
Remark 3.4.8 in \cite{la}, $\tilde{\rho}(F,\tilde{F})$ can also be
written as \ce \tilde{\rho}(F,\tilde{F})=\max_{i\neq
j}\|P_i\tilde{P}_j\|^{h/|\lambda_i-\lambda_j|}, \de where $P_i$
denotes the orthogonal projection onto $U_i$.

\subsection{L\'evy processes for $t\geq0$}\label{levypro1}

\bd\label{levy1} A process $L=(L_t)_{t\geq0}$ with $L_0=0$, a.s., is
a $d$-dimensional L\'evy process for $t\geq0$ if

(i) $L$ has independent increments; that is, $L_t-L_s$ is
independent of $L_v-L_u$ if $(u,v)\cap(s,t)=\emptyset$;

(ii) $L$ has stationary increments; that is, $L_t-L_s$ has the same
distribution as $L_v-L_u$ if $t-s=v-u>0$;

(iii) $L_t$ is right continuous with left limit. \ed

The characteristic function of $L_t$ is given by \ce
\mE\left(\exp\{i\<z,L_t\>\}\right)=\exp\{t\Psi(z)\}, \quad
z\in\mR^d. \de The function $\Psi: \mR^d\rightarrow\mathcal {C}$ is
called the characteristic exponent of the L\'evy process $L$. By the
L\'evy-Khintchine formula, there exist a nonnegative-definite
$d\times d$ matrix $Q$, a vector $\gamma\in\mR^d$, and a measure
$\nu$ on $\mR^d\setminus\{0\}$ satisfying \be
\int_{\mR^d\setminus\{0\}}(|u|^2\wedge1)\nu(\dif u)<\infty, \label{lemc} \ee such
that \be
\Psi(z)=-\frac{1}{2}\<z,Qz\>+i\<z,\gamma\>
+\int_{\mR^d\setminus\{0\}}\big(e^{i\<z,u\>}-1-i\<z,u\>1_{|u|\leq\delta}\big)\nu(\dif
u), \label{lkf} \ee where $\delta>0$ is a constant. Here $\nu$ is
called a L\'evy measure.

Set $\kappa_t:=L_t-L_{t-}$. Then $\kappa$ defines a stationary
$(\cF_0^t)_{t\geq0}$-adapted Poisson point process with values in
$\mR^d\setminus\{0\}$ and characteristic measure $\nu$ (c.f.
\cite{iw}). Let $N_{\kappa}((0,t],\dif u)$ be the counting measure
of $\kappa_{t}$, i.e., for $B\in\cB(\mR^d\setminus\{0\})$
$$
N_{\kappa}((0,t],B):=\#\{0<s\leq t: \kappa_s\in B\},
$$
where $\#$ denotes the cardinality of a set. The compensator measure
of $N_{\kappa}$ is given by
$$
\tilde{N}_{\kappa}((0,t],\dif u):=N_{\kappa}((0,t],\dif u)-t\nu(\dif
u).
$$
The L\'evy-It\^o theorem states that there exist  a vector
$b\in\mR^d$, a $d'$-dimensional $(\cF_0^t)_{t\geq0}$-Brownian motion
$W_t$, with $0\leq d'\leq d$, and  a $d\times d'$ matrix $A$, such
that $L$ can be represented as \be
L_t=bt+AW_t+\int_0^t\int_{|u|\leq\delta}u\tilde{N}_{\kappa}(\dif s, \dif u)
+\int_0^t\int_{|u|>\delta}uN_{\kappa}(\dif s, \dif u).
\label{lif1} \ee

\section{two motivated examples}\label{exam}

In the section, we will give two motivated examples.

\begin{example} \label{example1}
Consider a linear stochastic system in $\mR^2$ with L\'evy processes:
\be\left\{\begin{array}{l}
\dif X_t^1=\gamma X_t^1\dif t+X_t^1\dif L_t^1, \qquad X_0^1=x^1, \qquad t\geq0,\\
\dif X_t^2=-\gamma X_t^2\dif t+X_t^2\dif L_t^2, \quad X_0^2=x^2, ~\qquad
t\geq0,
\end{array}
\right.
\label{equa1}
\ee
where $\gamma>0$ is a constant, and $L^1_t, L^2_t$  are two
one-dimensional independent L\'evy processes with the same L\'evy
measure $\nu$ and the following L\'evy-It\^o representations,
\ce
L^1_t=\int_0^t\int_{|u|\leq\delta}u\tilde{N}_{\kappa^1}(\dif s, \dif u), \quad
L^2_t=\int_0^t\int_{|u|\leq\delta}u\tilde{N}_{\kappa^2}(\dif s, \dif u).
\de
Here, we require $0<\delta<1$. By the It\^o formula, we obtain the solution of Eq.(\ref{equa1})
\ce\left(\begin{array}{l}
X_t^1\\
X_t^2
\end{array}
\right) =\left(\begin{array}{l}
M_t^1 \quad 0\\
0 \quad M_t^2
\end{array}
\right) \left(\begin{array}{l}
x^1\\
x^2
\end{array}
\right), \de where \ce
&&M_t^1=\exp\bigg\{\Big[\gamma+\int_{|u|\leq\delta}\big(\log(1+u)-u\big)\nu(\dif
u)\Big]t+\int_0^t\int_{|u|\leq\delta}\log(1+u)
\tilde{N}_{\kappa^1}(\dif s, \dif u)\bigg\},\\
&&M_t^2=\exp\bigg\{\Big[-\gamma+\int_{|u|\leq\delta}\big(\log(1+u)-u\big)\nu(\dif
u)\Big]t+\int_0^t\int_{|u|\leq\delta}\log(1+u)
\tilde{N}_{\kappa^2}(\dif s, \dif u)\bigg\}. \de
Set
\ce
\varphi(t,\omega):=\left(\begin{array}{l}
M_t^1 \quad 0\\
0 \quad M_t^2
\end{array}
\right). \de
By the strong law of large numbers, it follows that \ce
\lim\limits_{t\rightarrow\infty}\frac{1}{t}\int_0^t\int_{|u|\leq\delta}\log(1+u)\tilde{N}_{\kappa^i}(\dif
s, \dif u)=0. \de Thus, we have \ce
\lim\limits_{t\rightarrow\infty}\big(\varphi(t,\omega)^*\varphi(t,\omega)\big)^{1/2t}
=\lim\limits_{t\rightarrow\infty}\left(\begin{array}{l}
(M^1_t)^2 \quad 0\\
0 \qquad (M^2_t)^2
\end{array}
\right)^{1/2t} =\left(\begin{array}{l}
M^1 \quad 0\\
0 \qquad M^2
\end{array}
\right)=:\Phi(\omega), \de where \ce
&&M^1=\exp\left\{\gamma+\int_{|u|\leq\delta}\big(\log(1+u)-u\big)\nu(\dif u)\right\},\\
&&M^2=\exp\left\{-\gamma+\int_{|u|\leq\delta}\big(\log(1+u)-u\big)\nu(\dif
u)\right\}. \de

It is obvious that $e^{\lambda_1}, e^{\lambda_2}$ are eigenvalues of $\Phi(\omega)$, where
\ce
&&\lambda_1=\gamma+\int_{|u|\leq\delta}\big(\log(1+u)-u\big)\nu(\dif u),\\
&&\lambda_2=-\gamma+\int_{|u|\leq\delta}\big(\log(1+u)-u\big)\nu(\dif u),
\de
and corresponding eigenspaces $U_1=\{(x^1,x^2)\in\mR^2| x^2=0\}$, $U_2=\{(x^1,x^2)\in\mR^2| x^1=0\}$.

Take $V_3=\{0\}$, $V_2=U_2$ and $V_1=\mR^2$. For each
$x\in\mR^2\setminus\{0\}$, \ce
\lambda(\omega,x)=\lambda_2~\Leftrightarrow~ x\in U_2=V_2\setminus V_3,\\
\lambda(\omega,x)=\lambda_1~\Leftrightarrow~ x\in U_1=V_1\setminus
V_2. \de

For $i=1, 2$, define a one-dimensional two-sided L\'evy process $\hat{L}_t^i$, for $t\in\mR$,
\ce
\hat{L}_t^i=\left\{\begin{array}{l}
L_t^i, ~\quad t\geq0,\\
\tilde{L}^{i}_t, \quad t<0,
\end{array}
\right.
\de
where $\tilde{L}^{i}_t$ is an independent copy of $-L_{(-t)-}^i$ and independent
of $L_t^j$ for $j=1, 2$ and $\tilde{L}_t^{k}$ for $k=1, 2, k\neq i$. Now consider the following system
\be\left\{\begin{array}{l}
\dif X_t^1=\gamma X_t^1\dif t+X_t^1\dif \hat{L}^1_t, \qquad X_0^1=x^1,\\
\dif X_t^2=-\gamma X_t^2\dif t+X_t^2\dif \hat{L}^2_t, \quad X_0^2=x^2,
\end{array}
\right. \label{equa1111} \ee where the corresponding stochastic
integrals are understood as the forward It\^o integrals for $t\geq0$
(\cite{iw}), and the backward It\^o integrals for $t\leq0$
(\cite{fk2}). Set
\ce
E_1:=U_1, \quad E_2:=U_2,
\de
and then by the similar deduction as that for $t\geq 0$ in the
first part of this example
\ce
\lim\limits_{t\rightarrow\pm\infty}\frac{1}{t}\log|\varphi(t,\omega)x|=\lambda_i(\omega)
\Longleftrightarrow x\in E_i(\omega)\setminus\{0\}, \quad i=1,2. \de
\end{example}

\begin{example} \label{example2}
Consider a linear stochastic system in $\mR^2$ with L\'evy processes:
\be\left\{\begin{array}{l}
\dif \bar{X}_t^1=\gamma \bar{X}_t^2\dif t+\bar{X}_t^1\dif L_t^1, \qquad \bar{X}_0^1=x^1, \qquad t\geq0,\\
\dif \bar{X}_t^2=-\gamma \bar{X}_t^1\dif t+\bar{X}_t^2\dif L_t^2, \quad \bar{X}_0^2=x^2, ~\qquad
t\geq0,
\end{array}
\right.
\label{equa2}
\ee
where $\gamma, L_t^1, L_t^2$ are the same to that in Eq.(\ref{equa1}). Although Eq.(\ref{equa2}) is a bit different from Eq.(\ref{equa1}) and is also linear, its solution and solution operator can not be explicitly expressed(\cite{p}). Thus, similar deduction to that in Example \ref{example1} is not done. And we only make use of another method to obtain its Lyapunov exponents.
\end{example}

\section{A MET for linear cocycles discontinuous in $t$}\label{met}

We first recall the following lemma for linear cocycles with discrete
time  (\cite[Theorem 3.4.11(A), p.153]{la}).

\bl({\bf MET for Linear Cocycle with Two-Sided Discrete Time})\label{met1}\\
Let \ce \varphi(n,\omega)=\left\{\begin{array}{c}
A(\theta^{n-1}\omega)\cdots A(\omega), ~\qquad\quad n\geq 1,\\
I, ~\qquad\qquad\qquad\qquad\quad\quad n=0,\\
A^{-1}(\theta^n\omega)\cdots A^{-1}(\theta^{-1}\omega), \quad\quad
n\leq -1,
\end{array}
\right. \de where $A: \Omega\mapsto Gl(d,\mR)$($Gl(d,\mR)$ denotes
the group of $d\times d$ invertible real matrices) is a strongly
measurable random invertible matrix and $\theta:
\Omega\mapsto\Omega$ is a measurable mapping with
$\theta^{-1}\Omega=\Omega$ and $\mP\theta^{-1}=\mP$. Assume \ce
\log^+\|A(\cdot)\|\in\cL^1(\Omega,\cF,\mP) ~\mbox{and}~
\log^+\|A^{-1}(\cdot)\|\in\cL^1(\Omega,\cF,\mP). \de Then there
exists an invariant set $\tilde{\Omega}$ of full measure such that
for $\omega\in\tilde{\Omega}$

(i) The limit
$\lim\limits_{n\rightarrow\infty}\big(\varphi(n,\omega)^*
\varphi(n,\omega)\big)^{1/2n}=:\Phi(\omega)\geq0$ exists.

(ii) Let
$e^{\lambda_{p(\omega)}(\omega)}<\cdots<e^{\lambda_1(\omega)}$ be
the different eigenvalues of $\Phi(\omega)$
($\lambda_{p(\omega)}>-\infty$) and let
$U_{p(\omega)}(\omega),\cdots,U_1(\omega)$ be the corresponding
eigenspaces with multiplicities $d_i(\omega):=\dim U_i(\omega)$.
Then \ce p(\theta\omega)=p(\omega), \quad
\lambda_i(\theta\omega)=\lambda_i(\omega), \quad
d_i(\theta\omega)=d_i(\omega), \de for $i=1,\dots,p(\omega)$.

(iii) Put $V_{p(\omega)+1}(\omega):=\{0\}$ and for
$i=1,\dots,p(\omega)$ \ce
V_i(\omega):=U_{p(\omega)}(\omega)\oplus\dots\oplus U_i(\omega), \de
so that \ce V_{p(\omega)}(\omega)\subset\dots\subset
V_i(\omega)\subset\dots\subset V_1(\omega)=\mR^d \de defines a
filtration of $\mR^d$. Then for each $x\in\mR^d\setminus\{0\}$ the
Lyapunov exponent \ce
\lambda(\omega,x):=\lim\limits_{n\rightarrow\infty}\frac{1}{n}\log|\varphi(n,\omega)x|
\de exists and \ce
\lambda(\omega,x)=\lambda_i(\omega)\Longleftrightarrow x\in
V_i(\omega)\setminus V_{i+1}(\omega), \de equivalently \ce
V_i(\omega)=\{x\in\mR^d: \lambda(\omega,x)\leq\lambda_i(\omega)\}.
\de

(iv) For all $x\in\mR^d\setminus\{0\}$ \ce
\lambda(\theta\omega,A(\omega)x)=\lambda(\omega,x), \de whence \ce
A(\omega)V_i(\omega)=V_i(\theta\omega) \de for
$i=1,\dots,p(\omega)$.

(v) If $(\Omega,\cF,\mP,\theta)$ is ergodic, $p(\omega)$ is a
constant on $\tilde{\Omega}$, and $\lambda_i(\omega)$ and
$d_i(\omega)$ are two constants on $\{\omega\in\Omega: p(\omega)\geq
i\}, ~i=1,\dots,d$.

(vi) For each $\omega\in\tilde{\Omega}$ there exists a splitting \ce
\mR^d=E_1(\omega)\oplus\dots\oplus E_{p(\omega)}(\omega) \de of
$\mR^d$ with $\dim E_i(\omega)=d_i(\omega)$ such that for
$i\in\{1,\dots,p(\omega)\}$,

(a) if $P_i(\omega): \mR^d\mapsto E_i(\omega)$ is the projection
onto $E_i(\omega)$ along $F_i(\omega):=\oplus_{j\neq i}E_j(\omega)$,
then \ce A(\omega)P_i(\omega)=P_i(\theta\omega)A(\omega), \de
equivalently \ce A(\omega)E_i(\omega)=E_i(\theta\omega), \de

(b) we have \ce
\lim\limits_{n\rightarrow\pm\infty}\frac{1}{n}\log|\varphi(n,\omega)x|=\lambda_i(\omega)
\Longleftrightarrow x\in E_i(\omega)\setminus\{0\}, \de

(c) convergence in (b) is uniform with respect to $x\in
E_i(\omega)\cap S$ for each fixed $\omega$, where $S=\{x\in\mR^d:
|x|=1\}$. \el

\medskip

Now we state and prove the following MET for linear cocycles discontinuous in $t$.

\bt({\bf MET for Two-Sided Linear Cocycle Discontinuous in $t$})\label{met2}\\
Let $\varphi: \mR\times\Omega\times\mR^d\mapsto\mR^d$ be a linear
cocycle discontinuous in $t$ over the metric DS
$(\Omega,\cF,\mP,(\theta_t)_{t\in\mR})$. Let $\varphi(t,\omega)\in
Gl(d,\mR)$. Assume that $\alpha^+\in\cL^1(\Omega,\cF,\mP)$ and
$\alpha^-\in\cL^1(\Omega,\cF,\mP)$, where \be
\alpha^+(\omega):=\sup\limits_{0\leq t\leq
1}\log^+\|\varphi(t,\omega)\|, \quad
\alpha^-(\omega):=\sup\limits_{0\leq t\leq
1}\log^+\|\varphi(t,\omega)^{-1}\|. \label{icv} \ee Then there
exists an invariant set $\tilde{\Omega}$ of full measure such that
for $\omega\in\tilde{\Omega}$ all statements of Lemma \ref{met1}
hold with $n, \theta$ and $A(\omega)$ replaced by $t, \theta_t$ and
$\varphi(t,\omega)$.
\et
\begin{proof}
{\bf Step 1. Measurability.} Because $\varphi(t,\omega)$ is
c\`adl\`ag for $t\geq0$, the random variables $\alpha^+$ and
$\alpha^-$ are $\cF$-measurable.

{\bf Step 2. Convergence of flags.} (i) For $t\in\mR_+$, there exist
two orthogonal matrices $G_t$ and $O_t$ such that
$$
\varphi(t)=G_tD_tO_t, \quad
D_t=diag\big(\delta_1(\varphi(t)),\cdots,\delta_d(\varphi(t))\big),
$$
where $\delta_i(\varphi(t))$ is the singular value of $\varphi(t)$
and $\delta_1(\varphi(t))
\geq\delta_2(\varphi(t))\geq\cdots\geq\delta_d(\varphi(t))>0$. By
Proposition 3.2.7 (iii) in \cite{la}, \be
\|\wedge^k\varphi(t)\|=\delta_1(\varphi(t))\cdots\delta_k(\varphi(t)),
\label{sdre} \ee where $\wedge^k\varphi(t)$ denotes the k-fold
exterior power of $\varphi(t)$.

Cocycle property for $\varphi(t)$ and Lemma 3.2.6 (v) in \cite{la}
allow us to get \ce
\wedge^k\varphi(t,\omega)=\big(\wedge^k\varphi(t-[t],\theta_{[t]}\omega)\big)\big(\wedge^k\varphi([t],\omega)\big)
\de and \ce
\wedge^k\varphi([t],\omega)&=&\big(\wedge^k\varphi(t-[t],\theta_{[t]}\omega)\big)^{-1}\big(\wedge^k\varphi(t,\omega)\big)\\
&=&\big(\wedge^k\varphi(t-[t],\theta_{[t]}\omega)^{-1}\big)\big(\wedge^k\varphi(t,\omega)\big).
\de Based on Proposition 3.2.7 (iii) in \cite{la}, it holds that \ce
\|\wedge^k\varphi(t,\omega)\|&\leq&\|\wedge^k\varphi([t],\omega)\|\|\varphi(t-[t],\theta_{[t]}\omega)\|^k\\
&\leq&\|\wedge^k\varphi([t],\omega)\|\left(\sup\limits_{0\leq s\leq
1}\|\varphi(s,\theta_{[t]}\omega)\|\right)^k \de and \ce
\|\wedge^k\varphi([t],\omega)\|&\leq&\|\varphi(t-[t],\theta_{[t]}\omega)^{-1}\|^k\|\wedge^k\varphi(t,\omega)\|\\
&\leq&\left(\sup\limits_{0\leq s\leq
1}\|\varphi(s,\theta_{[t]}\omega)^{-1}\|\right)^k\|\wedge^k\varphi(t,\omega)\|.
\de Thus, \ce
\frac{\log\|\wedge^k\varphi([t],\omega)\|}{t}-\frac{k\alpha^-(\theta_{[t]}\omega)}{t}&\leq&\frac{\log\|\wedge^k\varphi(t,\omega)\|}{t}\\
&\leq&\frac{\log\|\wedge^k\varphi([t],\omega)\|}{t}+\frac{k\alpha^+(\theta_{[t]}\omega)}{t},
\de where we have used the following two inequalities: \ce
&&\log\sup\limits_{0\leq s\leq
1}\|\varphi(s,\theta_{[t]}\omega)\|\leq\sup\limits_{0\leq s\leq
1}\log^+
\|\varphi(s,\theta_{[t]}\omega)\|,\\
&&\log\sup\limits_{0\leq s\leq
1}\|\varphi(s,\theta_{[t]}\omega)^{-1}\|\leq\sup\limits_{0\leq s\leq
1}\log^+ \|\varphi(s,\theta_{[t]}\omega)^{-1}\|. \de By (\ref{icv}),
we obtain that \ce
\lim\limits_{n\rightarrow\infty}\frac{\log\|\wedge^k\varphi(n,\omega)\|}{n}
=\lim\limits_{t\rightarrow\infty}\frac{\log\|\wedge^k\varphi(t,\omega)\|}{t}.
\de So, by Theorem 3.3.3(B) in \cite{la} for
$A(\omega)=\varphi(1,\omega)$ and $\theta=\theta_1$, there exist a
forward invariant set $\Omega_1\in\cF$ of full measure
($\Omega_1\subset\theta_1^{-1}\Omega_1$ and $\mP(\Omega_1)=1$) and
measurable functions $\gamma^{(k)}:\Omega\rightarrow\mR$, with
$(\gamma^{(k)})^+\in\cL^1(\Omega,\cF,\mP)$, such that \be
\lim\limits_{t\rightarrow\infty}\frac{\log\|\wedge^k\varphi(t,\omega)\|}{t}=\gamma^{(k)},
\qquad a.s.. \label{expl} \ee Combining (\ref{expl}) and
(\ref{sdre}), we have \ce
\lim\limits_{t\rightarrow\infty}D_t^{1/t}=diag(e^{\Lambda_1},\cdots,e^{\Lambda_d}),
\de where $\gamma^{(k)}=\Lambda_1+\cdots+\Lambda_k$.

Denote by $\lambda_1>\cdots>\lambda_p$ the distinct numbers among
the $\Lambda_i$. Let $d_i$ be the multiplicity of $\lambda_i$ for
$i=1,\dots,p$. Put \ce \Delta_i:=\lambda_i-\lambda_{i+1},
i=1,\dots,p-1, \quad \Delta:=\min\limits_{i=1,\dots,p-1}\Delta_i>0.
\de Let $U_i(t)$ be spanned by the group $\Sigma_i$ of those
eigenvectors of $(\varphi^*(t)\varphi(t))^{1/2t}=O^*_tD_t^{1/t}O_t$
corresponding to eigenvalues $\delta_{k(i)}(\varphi(t))^{1/t}$,
where
$\lim\limits_{t\rightarrow\infty}\delta_{k(i)}(\varphi(t))^{1/t}=e^{\lambda_i}$,
and \ce V_i(t):=U_p(t)\oplus\dots\oplus U_i(t), \quad i=1,\dots,p.
\de The sequence of subspaces of $\mR^d$ given by \ce
F(t)=\big(V_p(t),\dots,V_i(t),\dots,V_1(t)\big) \de forms a flag of
type \ce \tau=(d_p,d_p+d_{p-1},\cdots,d_p+\cdots+d_1=d). \de

Denote $h=\frac{\Delta}{d-1}$. So, by Section \ref{fm}, the distance
between $F(t)$ and $ F([t])$, in $F_\tau(d)$, is given by \ce
\tilde{\rho}\big(F(t), F([t])\big)=\max_{i\neq j\atop
i,j=1,\dots,p}\|P_i(t)P_j([t])\|^{h/|\lambda_i-\lambda_j|}, \de
where $P_i(t)$ denotes the orthogonal projection onto $U_i(t)$.

(ii) Next, we calculate $\tilde{\rho}(F(t),F([t]))$. If $i>j$,
$\lambda_i<\lambda_j$. Take a unit vector $x\in U_i([t])$ and
$y=P_j(t)x\in U_j(t)$. Thus, \ce
|\varphi(t,\omega)x|&=&|\varphi(t-[t],\theta_{[t]}\omega)\varphi([t],\omega)x|
\leq\|\varphi(t-[t],\theta_{[t]}\omega)\||\varphi([t],\omega)x|\\
&\leq&\|\varphi(t-[t],\theta_{[t]}\omega)\|\bar{\delta}_i(\varphi([t])),
\de and \ce
|\varphi(t,\omega)x|^2&=&|\varphi(t,\omega)y|^2+|\varphi(t,\omega)(x-y)|^2\\
&\geq&|\varphi(t,\omega)y|^2\geq\underline{\delta}_j(\varphi(t))^2|y|^2,
\de where \ce
\bar{\delta}_i(\varphi([t]))=\sup\limits_{\Sigma_i}\delta_{k(i)}(\varphi([t])),
\qquad
\underline{\delta}_j(\varphi(t))=\inf\limits_{\Sigma_j}\delta_{k(j)}(\varphi(t)),
\de with \ce
\lambda_i=\limsup\limits_{t\rightarrow\infty}\frac{1}{t}\log\bar{\delta}_i(\varphi([t])),
\qquad
\lambda_j=\limsup\limits_{t\rightarrow\infty}\frac{1}{t}\log\underline{\delta}_j(\varphi(t)).
\de Therefore, \ce
|y|=|P_j(t)P_i([t])x|\leq\frac{|\varphi(t,\omega)x|}{\underline{\delta}_j(\varphi(t))}
\leq\|\varphi(t-[t],\theta_{[t]}\omega)\|\frac{\bar{\delta}_i(\varphi([t]))}{\underline{\delta}_j(\varphi(t))},
\de and \ce
\|P_j(t)P_i([t])\|&\leq&\|\varphi(t-[t],\theta_{[t]}\omega)\|\frac{\bar{\delta}_i(\varphi([t]))}{\underline{\delta}_j(\varphi(t))}\\
&\leq&\sup\limits_{0\leq
s\leq1}\|\varphi(s,\theta_{[t]}\omega)\|\frac{\bar{\delta}_i(\varphi([t]))}{\underline{\delta}_j(\varphi(t))}.
\de Moreover, \ce
\limsup\limits_{t\rightarrow\infty}\frac{1}{t}\log\|P_j(t)P_i([t])\|
\leq\limsup\limits_{t\rightarrow\infty}\frac{1}{t}\log\sup\limits_{0\leq
s\leq1}\|\varphi(s,\theta_{[t]}\omega)\|-|\lambda_i-\lambda_j|. \de
Since $\log\sup\limits_{0\leq
s\leq1}\|\varphi(s,\omega)\|\leq\sup\limits_{0\leq s\leq
1}\log^+\|\varphi(s,\omega)\|\in\cL^1(\Omega,\cF,\mP)$, we deduce
that \ce
\limsup\limits_{t\rightarrow\infty}\frac{1}{t}\log\sup\limits_{0\leq
s\leq1}\|\varphi(s,\theta_{[t]}\omega)\|\leq0. \de Thus, \be
\limsup\limits_{t\rightarrow\infty}\frac{1}{t}\log\|P_j(t)P_i([t])\|\leq-|\lambda_i-\lambda_j|.
\label{any1} \ee

If $i<j$, $\lambda_i>\lambda_j$. By the same deduction as above, we
obtain \ce
\|P_i([t])P_j(t)\|\leq\|\varphi(t-[t],\theta_{[t]}\omega)^{-1}\|\frac{\bar{\delta}_j(\varphi(t))}
{\underline{\delta}_i(\varphi([t]))} \de and \be
\limsup\limits_{t\rightarrow\infty}\frac{1}{t}\log\|P_j(t)P_i([t])\|\leq-|\lambda_i-\lambda_j|.
\label{any2} \ee

Combining (\ref{any2}) and (\ref{any1}), we get \be
\limsup_{t\rightarrow\infty}\frac{1}{t}\log\tilde{\rho}\left(F(t),F([t])\right)\leq-h.
\label{e1} \ee

(iii) By Lemma 3.4.9 in \cite{la}, there exists a flag
$F=\big(V_p,\dots,V_i,\dots,V_1\big)$ of type $\tau$ such that \be
\limsup_{n\rightarrow\infty}\frac{1}{n}\log\tilde{\rho}\left(F(n),F\right)\leq-h.
\label{e2} \ee By (\ref{e1}) and (\ref{e2}), we have \ce
\limsup_{t\rightarrow\infty}\frac{1}{t}\log\tilde{\rho}\left(F(t),F\right)
&\leq&\limsup_{t\rightarrow\infty}\frac{1}{t}\log\Big(\tilde{\rho}\big(F(t),F([t])\big)+\tilde{\rho}\big(F([t]),F\big)\Big)\\
&\leq&-h. \de

{\bf Step 3. Lyapunov exponent.} If $t=n+s$ with $s\in(0,1),
n\in\mN$, then \ce
\|\varphi(s,\theta_n\omega)\||\varphi(n,\omega)x|\geq|\varphi(t,\omega)x|
\geq\|\varphi(s,\theta_n\omega)^{-1}\|^{-1}|\varphi(n,\omega)x|, \de
and therefore \ce \alpha^+(\theta_n\omega)+\log|\varphi(n,\omega)x|
\geq\log|\varphi(t,\omega)x|\geq\log|\varphi(n,\omega)x|-\alpha^-(\theta_n\omega).
\de Since
$\lim\limits_{n\rightarrow\infty}n^{-1}\alpha^+(\theta_n\omega)
=\lim\limits_{n\rightarrow\infty}n^{-1}\alpha^-(\theta_n\omega)=0$
with probability 1, one has \ce
\lim\limits_{t\rightarrow\infty}\frac{1}{t}\log|\varphi(t,\omega)x|=
\lim\limits_{n\rightarrow\infty}\frac{1}{n}\log|\varphi(n,\omega)x|.
\de By Lemma \ref{met1} (iii), the statement in (iii) holds.

{\bf Step 4. Invariancy.} For   $x\in\mR^d\setminus\{0\}$, \ce
\lambda(\theta_t\omega,\varphi(t,\omega)x)
&=&\limsup_{s\rightarrow\infty}\frac{1}{s}\log|\varphi(s,\theta_t\omega)\varphi(t,\omega)x|\\
&=&\limsup_{s\rightarrow\infty}\frac{1}{s}\log|\varphi(s+t,\omega)x|\\
&=&\limsup_{s\rightarrow\infty}\frac{1}{s+t}\log|\varphi(s+t,\omega)x|\cdot\frac{s+t}{s}\\
&=&\lambda(\omega,x). \de

{\bf Step 5. The flag for negative time.} For $t\in\mR_-$, cocycle
property for $\varphi(t,\omega)$ infers that \ce
\varphi(t,\omega)=\varphi(-t,\theta_t\omega)^{-1}. \de Let
$\delta_1\big(\varphi(-t,\theta_t\omega)\big)\geq\delta_2\big(\varphi(-t,\theta_t\omega)\big)\geq\cdots\geq
\delta_d\big(\varphi(-t,\theta_t\omega)\big)>0$ be singular values
of $\varphi(-t,\theta_t\omega)$, and then these singular values
$\delta_1\big(\varphi(t)\big)\geq\delta_2\big(\varphi(t)\big)\geq\cdots\geq
\delta_d\big(\varphi(t)\big)$ of $\varphi(t)$ satisfy
$$
\delta_k\big(\varphi(t)\big)=\delta_{d+1-k}\big(\varphi(-t,\theta_t\omega)\big)^{-1},
$$
for $k=1,2,\cdots,d$. By the same deduction as that in Step 2, we
get that \ce
\lim\limits_{n\rightarrow-\infty}\frac{\log\|\wedge^k\varphi(n,\omega)\|}{|n|}
=\lim\limits_{t\rightarrow-\infty}\frac{\log\|\wedge^k\varphi(t,\omega)\|}{|t|}.
\de So, by Theorem 3.3.10(A) in \cite{la} for
$A^{-1}(\theta^{-1}\omega)=\varphi(-1,\omega)$ and
$\theta^{-1}=\theta_{-1}$, on a invariant set $\Omega_2\in\cF$ of
full measure \be
\lim\limits_{n\rightarrow-\infty}\frac{\log\|\wedge^k\varphi(n,\omega)\|}{|n|}
=\gamma^{(d-k)}-\gamma^{(d)}, \qquad a.s.. \label{expl1} \ee
Combining (\ref{expl1}) and (\ref{sdre}), we have \ce
\lim\limits_{t\rightarrow-\infty}D_t^{1/t}=diag(e^{\Lambda^-_1},\cdots,e^{\Lambda^-_d}),
\de where $\Lambda^-_k=-\Lambda_{d+1-k}$.

Denote by $\lambda^-_1>\cdots>\lambda^-_p$ the distinct  numbers
among   $\Lambda^-_i$. Let $d^-_i$ be the multiplicity of
$\lambda^-_i$ for $i=1,\dots,p$. Then \ce
\lambda^-_k=-\lambda_{p+1-k}, \qquad d^-_k=d_{p+1-k}. \de By the
same deduction as in Step 2, we get a flag
$F^-=\big(V^-_p,V^-_{p-1},\cdots, V^-_1\big)$ of type
$\tau^-=(d^-_p,d^-_p+d^-_{p-1},\cdots,d^-_p+\cdots+d^-_1=d).$

{\bf Step 6. Oseledets spaces.} Let \ce E_i=V_i\cap V^-_{p+1-i},
\qquad i=1,2,\cdots,p. \de So, by the proof of Theorem 3.4.11(A) in
\cite{la}, $E_1, E_2,\cdots, E_p$, which form a
  splitting of $\mR^d$,   are Oseledets spaces.

In the following, we examine the properties of Oseledets spaces
$E_i$.

(i) For $t\geq0$ by Step 4, and for $t\leq0$, similar to that in
Step 4, we obtain \ce
\varphi(t,\omega)E_i&=&\varphi(t,\omega)V_i\cap\varphi(t,\omega)V^-_{p+1-i}\\
&=&V_i(\theta_t\omega)\cap V^-_{p(\theta_t\omega)+1-i}(\theta_t\omega)\\
&=&E_i(\theta_t\omega). \de

(ii) For $t\geq0$ by Step 3, and for $t\leq0$   similar to that in
Step 3, it holds that \ce
\lim\limits_{t\rightarrow\pm\infty}\frac{1}{t}\log|\varphi(t,\omega)x|=
\lim\limits_{n\rightarrow\pm\infty}\frac{1}{n}\log|\varphi(n,\omega)x|.
\de Thus, by Lemma \ref{met1}(vi)(b), we have \ce
\lim\limits_{t\rightarrow\pm\infty}\frac{1}{t}\log|\varphi(t,\omega)x|=\lambda_i(\omega)
\Longleftrightarrow x\in E_i(\omega)\setminus\{0\}. \de The proof is
thus completed.
\end{proof}

\medskip

Next, we apply the above theorem to the two examples in Section \ref{exam}.

\begin{example} \label{example3}
(continuity of Example \ref{example1}) Set
\ce
\varphi(t,\omega):=\left(\begin{array}{l}
M_t^1 \quad 0\\
0 \quad M_t^2
\end{array}
\right),
\de
and then $\varphi(t,\omega)\in Gl(2,\mR)$. Moreover, by the properties of L\'evy processes it is easy to justify
that $\varphi(t,\omega)$ is a linear cocycle discontinuous in $t$. Thus, we only need to prove that $\varphi(t,\omega)$ satisfies
two integrability conditions in order to use Theorem \ref{met2}. Rewrite Eq.(\ref{equa1}) as
\ce\left\{\begin{array}{l}
\dif X_t=aX_t\dif t+\sigma^1X_t\dif L_t^1+\sigma^2X_t\dif L_t^2,\\
X_0=x,
\end{array}
\right. \de where \ce &&X_t=\left(\begin{array}{l}
X_t^1\\
X_t^2
\end{array}
\right), x=\left(\begin{array}{l}
x^1\\
x^2
\end{array}
\right),  a=\left(\begin{array}{l}
\gamma \qquad 0\\
0 \quad -\gamma
\end{array}
\right), \sigma^1=\left(\begin{array}{l}
1 \quad 0\\
0 \quad 0
\end{array}
\right), \sigma^2=\left(\begin{array}{l}
0 \quad 0\\
0 \quad 1
\end{array}
\right). \de Applying the It\^o formula to $\log|X_t|$, we infer
that \ce
\log|X_t|&=&\log|X_0|+\int_0^t\frac{X_s^ia_{ik}X_s^k}{|X_s|^2}\dif s+\int_0^t\int_{|u|\leq\delta}\left[\log|X_{s-}+u\sigma^jX_{s-}|-\log|X_{s-}|\right]\tilde{N}_{\kappa^j}(\dif s, \dif u)\\
&&+\int_0^t\int_{|u|\leq\delta}\left[\log|X_{s-}+u\sigma^jX_{s-}|-\log|X_{s-}|-u\frac{X_{s-}^i\sigma^j_{ik}X_{s-}^k}{|X_{s-}|^2}\right]\nu(\dif
u)\dif s. \de So, \ce \sup\limits_{0\leq t\leq
1}\log^+\|\varphi(t,\omega)\|\leq\sup\limits_{|x|=1}\sup\limits_{0\leq
t\leq 1}\big|\log|X_t|\big| \leq \gamma+I_1+I_2, \de where \ce
&&I_1=\sup\limits_{|x|=1}\sup\limits_{0\leq t\leq
1}\left|\int_0^t\int_{|u|\leq\delta}\left[\log|X_{s-}+u\sigma^jX_{s-}|-\log|X_{s-}|\right]
\tilde{N}_{\kappa^j}(\dif s, \dif u)\right|,\\
&&I_2=\sup\limits_{|x|=1}\sup\limits_{0\leq t\leq 1}\bigg|\int_0^t\int_{|u|\leq\delta}\Big[\log|X_{s-}+u\sigma^jX_{s-}|-\log|X_{s-}|
-u\frac{X_{s-}^i\sigma^j_{ik}X_{s-}^k}{|X_{s-}|^2}\Big]\nu(\dif
u)\dif s\bigg|.
\de
For $I_1$, by BDG inequality, mean value theorem and H\"older's inequality, we have
\ce
\mE I_1&\leq&\sum\limits_{j=1}^2\mE\left[\sup\limits_{|x|=1}\int_0^1\int_{|u|\leq\delta}(\log|X_{s-}+u\sigma^jX_{s-}|-\log|X_{s-}|)^2
N_{\kappa^j}(\dif s, \dif u)\right]^{\frac{1}{2}}\\
&\leq&\sum\limits_{j=1}^2\mE\left[\sup\limits_{|x|=1}\int_0^1\int_{|u|\leq\delta}\frac{|u|^2}{(1-|u|)^2}
N_{\kappa^j}(\dif s, \dif u)\right]^{\frac{1}{2}}\\
&\leq&\sum\limits_{j=1}^2\left[\mE\left(\int_0^1\int_{|u|\leq\delta}\frac{|u|^2}{(1-\delta)^2}
N_{\kappa^j}(\dif s, \dif u)\right)\right]^{\frac{1}{2}}\\
&=&\frac{2}{1-\delta}\left[\int_{|u|\leq\delta}|u|^2\nu(\dif u)\right]^{\frac{1}{2}}.
\de
For $I_2$, by mean value theorem, it holds that
\ce
\mE I_2&\leq&\sum\limits_{j=1}^2\mE\bigg[\sup\limits_{|x|=1}\sup\limits_{0\leq t\leq 1}\int_0^t\int_{|u|\leq\delta}\Big|\log|X_{s-}+u\sigma^jX_{s-}|-\log|X_{s-}|
-u\frac{X_{s-}^i\sigma^j_{ik}X_{s-}^k}{|X_{s-}|^2}\Big|\nu(\dif u)\dif s\bigg]\\
&\leq&\sum\limits_{j=1}^2\mE\bigg[\sup\limits_{|x|=1}\sup\limits_{0\leq t\leq 1}3\int_0^t\int_{|u|\leq\delta}\frac{|u|^2}{(1-|u|)^2}\nu(\dif u)\dif s\bigg]
\leq\frac{6}{(1-\delta)^2}\int_{|u|\leq\delta}|u|^2\nu(\dif u).
\de Thus, \ce \mE\alpha^+=\mE\left(\sup\limits_{0\leq t\leq
1}\log^+\|\varphi(t,\omega)\|\right)<\infty. \de

As in linear algebra, we find the inverse \ce
\varphi(t,\omega)^{-1}=\left(\begin{array}{l}
(M_t^1)^{-1} \quad 0\\
0 \qquad\quad (M_t^2)^{-1}
\end{array}
\right). \de Simple calculations lead to \ce
\log\|\varphi(t,\omega)^{-1}\|&=&\frac{1}{2}\log\big((M_t^1)^{-2}+(M_t^2)^{-2}\big)
=\frac{1}{2}\log\big((M_t^1)^2+(M_t^2)^2\big)-\log M_t^1-\log M_t^2\\
&=&\log\|\varphi(t,\omega)\|-\log M_t^1-\log M_t^2
\leq\log\|\varphi(t,\omega)\|+|\log M_t^1|+|\log M_t^2|. \de By
Jensen's inequality, we obtain \ce \mE\left(\sup\limits_{0\leq t\leq
1}\log^+\|\varphi(t,\omega)^{-1}\|\right)
&\leq&\mE\left(\sup\limits_{0\leq t\leq 1}\log^+\|\varphi(t,\omega)\|\right)+\mE\left(\sup\limits_{0\leq t\leq 1}|\log M_t^1|\right)\\
&&+\mE\left(\sup\limits_{0\leq t\leq 1}|\log M_t^2|\right).
\de
For the second term in the right hand side of the
above inequality, it follows from BDG inequality and the H\"older
inequality that \ce \mE\left(\sup\limits_{0\leq t\leq 1}|\log
M_t^1|\right)
&\leq&\mE\left(\sup\limits_{0\leq t\leq 1}\Big|\gamma+\int_{|u|\leq\delta}\big(\log(1+u)-u\big)\nu(\dif u)\Big|t\right)\\
&&+\mE\left(\sup\limits_{0\leq t\leq
1}\left|\int_0^t\int_{|u|\leq\delta}\log(1+u)
\tilde{N}_{\kappa^1}(\dif s, \dif u)\right|\right)\\
&\leq&\left|\gamma+\int_{|u|\leq\delta}\big(\log(1+u)-u\big)\nu(\dif u)\right|\\
&&+\mE\left(\int_0^1\int_{|u|\leq\delta}(\log(1+u))^2
N_{\kappa^1}(\dif s, \dif u)\right)^\frac{1}{2}\\
&\leq&\gamma+\int_{|u|\leq\delta}\left|\log(1+u)-u\right|\nu(\dif u)\\
&&+\left(\mE\left(\int_0^1\int_{|u|\leq\delta}(\log(1+u))^2
N_{\kappa^1}(\dif s, \dif u)\right)\right)^\frac{1}{2}\\
&=&\gamma+\int_{|u|\leq\delta}\left|\log(1+u)-u\right|\nu(\dif u)\\
&&+\left(\int_{|u|\leq\delta}\big(\log(1+u)\big)^2\nu(\dif
u)\right)^\frac{1}{2}. \de Since $|\log(1+u)-u|\leq C|u|^2$ and
$\log^2(1+u)\leq C|u|^2$ for $|u|\leq\delta$, we thus have by
(\ref{lemc})
$\mE\left(\sup\limits_{0\leq t\leq 1}|\log M_t^1|\right)<\infty.$
Similarly, we also   obtain $\mE\left(\sup\limits_{0\leq t\leq
1}|\log M_t^2|\right)<\infty.$ Thus, \ce
\mE\alpha^-=\mE\left(\sup\limits_{0\leq t\leq
1}\log^+\|\varphi(t,\omega)^{-1}\|\right)<\infty.
\de

Finally, by Theorem \ref{met2} the linear structure of $\Phi(\omega)$ is given.
\end{example}

\begin{example} \label{example4}
(continuity of Example \ref{example2}) Rewrite Eq.(\ref{equa2}) as
\ce\left\{\begin{array}{l}
\dif \bar{X}_t=\bar{a}\bar{X}_t\dif t+\sigma^1\bar{X}_t\dif L_t^1+\sigma^2\bar{X}_t\dif L_t^2,\\
\bar{X}_0=x,
\end{array}
\right. \de where \ce &&\bar{X}_t=\left(\begin{array}{l}
\bar{X}_t^1\\
\bar{X}_t^2
\end{array}
\right), \bar{x}=\left(\begin{array}{l}
\bar{x}^1\\
\bar{x}^2
\end{array}
\right),  \bar{a}=\left(\begin{array}{l}
0 \qquad \gamma\\
-\gamma \quad 0
\end{array}
\right),
\de
and $\sigma^1, \sigma^2$ are the same to that in Example \ref{example3}. By \cite[Theorem 6, Page 249]{p}, the
equation has a unique solution denoted by $\bar{X}_t(\bar{x})$. Set
\ce
\bar{\varphi}(t,\omega)\bar{x}:=\bar{X}_t(\bar{x}).
\de
And then by expanding Eq.(\ref{equa2}) from $t\geq0$ to $t\leq0$ as that in Example \ref{example1}, we obtain
$\bar{\varphi}(t,\omega)$ for $t\leq0$. Define
\ce
\bar{\varphi}(t,\omega)^{-1}:=\bar{\varphi}(-t,\theta_t\omega), \quad t\geq0.
\de
Thus $\bar{\varphi}(t,\omega)\in Gl(2,\mR)$. Moreover, it follows from the properties of L\'evy processes
that $\bar{\varphi}(t,\omega)$ is a linear cocycle discontinuous in $t$. Next, we justify that $\bar{\varphi}(t,\omega)$ satisfies
two integrability conditions. Applying the It\^o formula to $\log|X_t|$, one could obtain that
\ce
\log|\bar{X}_t|&=&\log|\bar{X}_0|+\int_0^t\frac{\bar{X}_s^i\bar{a}_{ik}\bar{X}_s^k}{|\bar{X}_s|^2}\dif s+\int_0^t\int_{|u|\leq\delta}\left[\log|\bar{X}_{s-}+u\sigma^j\bar{X}_{s-}|-\log|\bar{X}_{s-}|\right]\tilde{N}_{\kappa^j}(\dif s, \dif u)\\
&&+\int_0^t\int_{|u|\leq\delta}\left[\log|\bar{X}_{s-}+u\sigma^j\bar{X}_{s-}|-\log|\bar{X}_{s-}|-u\frac{\bar{X}_{s-}^i\sigma^j_{ik}\bar{X}_{s-}^k}{|\bar{X}_{s-}|^2}\right]\nu(\dif
u)\dif s\\
&=&\log|\bar{X}_0|+\int_0^t\int_{|u|\leq\delta}\left[\log|\bar{X}_{s-}+u\sigma^j\bar{X}_{s-}|-\log|\bar{X}_{s-}|\right]\tilde{N}_{\kappa^j}(\dif s, \dif u)\\
&&+\int_0^t\int_{|u|\leq\delta}\left[\log|\bar{X}_{s-}+u\sigma^j\bar{X}_{s-}|-\log|\bar{X}_{s-}|-u\frac{\bar{X}_{s-}^i\sigma^j_{ik}\bar{X}_{s-}^k}{|\bar{X}_{s-}|^2}\right]\nu(\dif
u)\dif s.
\de
By the same deduction to that in Example \ref{example3}, we have that
\ce
\mE\alpha^+=\mE\left(\sup\limits_{0\leq t\leq1}\log^+\|\bar{\varphi}(t,\omega)\|\right)<\infty.
\de
The fact $\bar{\varphi}(t,\omega)^{-1}=\bar{\varphi}(-t,\theta_t\omega)$ and similar estimate to above admit us to get that
\ce
\mE\alpha^-=\mE\left(\sup\limits_{0\leq t\leq
1}\log^+\|\bar{\varphi}(t,\omega)^{-1}\|\right)<\infty.
\de
Thus, by Theorem \ref{met2}, $\bar{\varphi}(t,\omega)$ is exponentially stable.
\end{example}

\end{document}